\def\IS{{\Bbb S}}

\def\IZ{{\Bbb Z}} 
 
\def\IC{\Bbb C} 
\def\ID{{\Bbb D}}
\def\tr{{\rm tr}} 
\def\oC{\hat{\IC}}

\documentclass[12pt]{article} 

\usepackage[psamsfonts]{amssymb} 
\usepackage{amsfonts,amsmath}

\usepackage{epsfig,multicol}

\newtheorem{theorem}{Theorem}
\newtheorem{lemma}{Lemma}
\newtheorem{corollary}{Corollary}

\title{Nondiscrete parabolic characters of the free group $F_2$:  supergroup density and Nielsen classes in the complement of the Riley slice.}
\author{ Gaven J. Martin \thanks {Research supported in part by the New Zealand Marsden Fund. \newline \newline Keywords:  Kleinian group,  two bridge, knot complement, parabolic character \newline 2010 MSC 30F40; 32G05 } } 
\date{} 
\begin{document}

\maketitle

\begin{abstract} A parabolic representation of the free group $F_2$ is one in which the images of both generators are parabolic elements of $PSL(2,\IC)$. 
The Riley slice is a closed subset ${\cal R}\subset \IC$ which is a model  for the parabolic, discrete and faithful characters of $F_2$.  The complement of the Riley slice is a bounded Jordan domain within which there are isolated points,  accumulating only at the boundary,  corresponding to parabolic discrete and faithful representations of rigid subgroups of $PSL(2,\IC)$.  Recent work of  Aimi, Akiyoshi, Lee, Oshika, Parker, Lee, Sakai, Sakuma \& Yoshida, have topologically identified all these groups.  Here we give the first identified substantive properties of the nondiscrete representations and prove a supergroup density theorem: given any irreducible parabolic representation $\rho_*:F_2\to PSL(2,\IC)$ whatsoever,   any non-discrete parabolic representation $\rho_0$ has an arbitrarily small perturbation $\rho_\epsilon$ so that $\rho_\epsilon(F_2)$ contains a conjugate of $\rho_*(F_2)$ as a proper subgroup. This implies that if $\Gamma_*$ is any nonelementary group generated by two parabolic elements (discrete or otherwise) and $\gamma_0$ is any point in the complement of the Riley slice,  then in any neighbourhood of $\gamma$ there is a point corresponding to a nonelementary group generated by two parabolics with a conjugate of $\Gamma_*$ as a proper subgroup.  Using these ideas we then show that there are nondiscrete parabolic representations with an arbitrarily large number of distinct Nielsen classes of parabolic generators.   \end{abstract}

\section{Introduction.}  

 Given a finitely generated group $\Gamma$,
 \[ \Gamma=\langle g_1,g_2,\ldots,g_n| r_i=1,i=1,\ldots m \rangle,\]   the set of all representations of $\Gamma$ in  $PSL(2,\IC)$ is denoted
by $R(\Gamma)$, and it is called the {\em variety of representations} of $\Gamma$.  $R(\Gamma)$ has a natural structure as an affine algebraic set over the complex numbers $\IC$ via the natural embedding 
\[ R(\Gamma)\ni  \rho\mapsto \big(\rho(g_1),\ldots,\rho(g_n)\big) \in PSL(2,\IC)\times \cdots \times PSL(2,\IC)  \] 
with the defining equation induced by the relations of $\Gamma=\langle g_1,g_2,\ldots,g_n| r_i=1,i=1,\ldots m \rangle$.  Notice that by Scott's Theorem a finitely generated subgroup of $PSL(2,\IC)$ is finitely presented  \cite{Scott}, so there are finitely many defining equations.  The action of $PSL(2,\IC)$ by conjugation on $R(\Gamma)$ is algebraic,  but the quotient $R(\Gamma)/PSL(2,\IC)$ may not be Hausdorff and so one typically considers the algebraic quotient of invariant theory since the group $PSL(2,\IC)$ is reductive.  This quotient is denoted $X(\Gamma)$ -- an affine algebraic set with a regular map inducing an isomorphism $R(\Gamma)\to X(\Gamma)$.  Given a representation $\rho\in R(\Gamma)$ its {\em character} is the map
\[ \chi_\rho:\Gamma\to \IC, \hskip15pt  \Gamma\ni g \to \chi_\rho(g) =\tr^2(g) \in \IC \] 
 In what follows we denote the free group on two generators (typically $a$ and $b$) by $F_2$,  
 \[ F_2=\langle a,b \rangle \]
 and it is well known that $X(F_2)\cong \IC^3$.  In fact,  given a representation $\rho:F_2\to PSL(2,\IC)$ the three complex numbers
 \[ \gamma_\rho = \tr [\rho(a),\rho(b)] -2,\;\; \beta_\rho(a)=\tr^2(\rho(a))-4\;\;{\rm and} \;\;  \beta_\rho(b)=\tr^2(\rho(b))-4 \]
 uniquely determine $\rho(F_2)\subset PSL(2,\IC)$ up to conjugacy away from the locus $\{\gamma_\rho=0\}$, (see e.g. \cite{GM1,MM}, but really a consequence of the Fricke identities \cite{FK}). It is an elementary fact, (see e.g. \cite{B}) that $\gamma_\rho=0$ if and only if the representation $\rho$ is reducible.  
 
 \medskip
 
 We call the triple $(\gamma_\rho,\beta_\rho(a),\beta_\rho(b))$ the {\em principal character} of $\rho$ and a {\em parabolic representation} of $F_2$ is $\rho\in{\cal R}(F_2)$ for which $\beta_\rho(a)=\beta_\rho(b)=0$. 

\begin{theorem} An irreducible parabolic representation $\rho\in {\cal R}(F_2)$ is uniquely determined up to conjugacy by $\gamma_\rho \in\IC\setminus \{0\}$ and the image is conjugate to the group $\langle f,g \rangle$ where
\begin{equation}\label{fg} f= \pm \left[\begin{array}{cc} 1 & 1 \\ 0& 1 \end{array}\right], \;\; g = \pm \left[\begin{array}{cc} 1 & 0 \\ z_\rho & 1 \end{array}\right], \;\; {\rm and}\;\;  \gamma_\rho=z_\rho^2. \end{equation}
\end{theorem}
Notice that $\gamma_\rho$ is an invariant of the Nielsen class of generating pairs.  In this way the set of conjugacy classes irreducible parabolic representations in $X(F_2)$ (call it $P(F_2)\subset X(F_2)$) is naturally identified with $\IC\setminus \{0\}$,  and this identification gives a Hausdorff topology on $P(F_2)$,  via the identification $P(F_2) \ni \rho \leftrightarrow \gamma_\rho \in \IC\setminus\{0\} $.

\medskip

So identified,  the space $P(F_2)$ is further partitioned into three sets corresponding to the discrete and free representations,  discrete but not free,  and the nondiscrete representations.  We discuss the first two briefly as our main results pertain to the third class.  

\subsection{Discrete and free points of $P(F_2)$.}

We need to define the Riley slice.  The group $PSL(2,\IC)$ acts on the $2$-sphere by linear fractional transformations - that is as the group of conformal homeomorphisms of $\IS^2$ - M\"obius group, $M(\IS^2)$,
\[ PSL(2,\IC) \ni \pm \left[\begin{array}{cc} a & b \\ c &  d \end{array}\right] \leftrightarrow \frac{az+b}{cz+d} \in  M(\IS^2). \]
Given a discrete subgroup $\Gamma\subset M(\IS^2)$,  the action of $\Gamma$ on $\IS^2$ partitions the sphere into two regions.  The {\em ordinary set} $\Omega_\Gamma$,  an open set where the action is discontinuous,  and its complement the {\em limit set} $\Lambda_\Gamma$ where the action is chaotic.  It is possible that $\Omega_\Gamma=\emptyset$.  However when $\Omega_\Gamma\neq \emptyset$,  $\Omega_\Gamma/\Gamma=\Sigma$ is a Riemann surface. 
Define
\begin{equation}\label{Gz} \Gamma_z=\Big\langle \left[\begin{array}{cc} 1 & 1 \\ 0& 1 \end{array}\right],\left[\begin{array}{cc} 1 & 0 \\ z  & 1 \end{array}\right] \Big\rangle \subset PSL(2,\IC)\end{equation}
and then the Riley slice
\[ {\cal R} = \{z\in\IC :  \mbox{$\Omega_{\Gamma_z}/\Gamma_z$ is the four times punctured sphere} \}.\]
We know that all the groups $\Gamma_z$,  $z\in {\cal R}$,  are quasiconformally conjugate (see, for example, \cite{Bers}) so ${\cal R}$ is connected and is in fact the quotient of the Teichm\"uller space of the four times punctured sphere by the subgroup of the mapping class group generated by Dehn twists about a curve which separates the punctures in pairs, \cite{Kra}. This Teichm\"uller space is simply connected and so ${\cal R}$ is topologically an annulus.  At each point $\Lambda(\Gamma_z)$ is a Cantor set and $\Gamma_z$ is isomorphic to $F_2$, the free group on two generators. J\o rgensen's theorem on algebraic limits \cite{Jorgensen} shows that ${\overline {\cal R}}$ is precisely the set of discrete parabolic representations of $F_2$. Indeed, if $z\in \partial{\cal R}$,  then $\Omega_{\Gamma_z}/\Gamma_z$ is a pair of triply punctured spheres are called `cusps' - the cusps of Figure 1 below \cite{KS,KoS},  or if  $z$ is not a cusp,  then $\Omega_z=\emptyset$ and $\Gamma_z$  is called degenerate.  The Riley slice is symmetric under complex conjugation and as $\gamma(f,g)=z^2$ we have $-{\cal R}={\cal R}$.  

Lyndon and Ullman were the first to make substantial advances in the study of this set, \cite{LU}, though there was earlier work of Sanov (1947), Chang, Jennings \& Ree (1958) and Brenner (1961).  But it was Riley's computational investigations and the remarkable connections to knot theory which really aroused mathematical interest. The important results of  Ohshika and Miyachi \cite{OM}, Keen and Series, \cite{KS} and also of Akiyoshi, Sakuma, Wada \& Yamashita \cite{ASWY} have shown us that the interior, $int({\cal R})$,  is conformally equivalent to $\IC\setminus \ID \cong \ID\setminus \{0\}$ and $\partial {\cal R}$ is a topological circle.  Each point of $int({\cal R})$ corresponds to a geometrically finite group $\Gamma$ generated by two parabolics. 

\scalebox{0.75}{\includegraphics[viewport= -5 560 580 800]{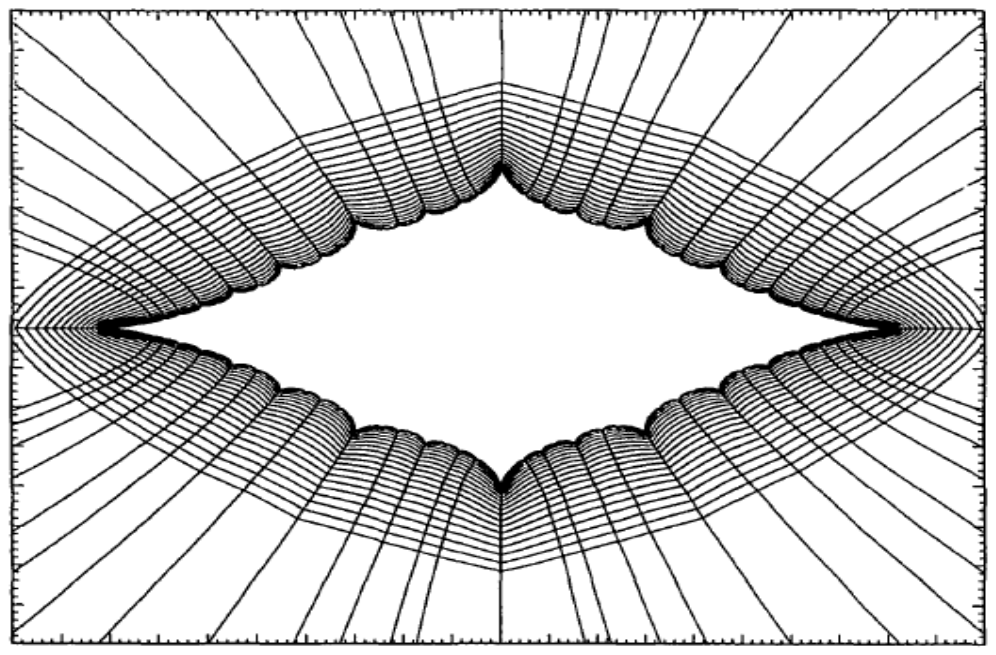}}\\
{\bf Figure 1.}  {\em The Riley slice (shaded region) is symmetric about the origin and its complement meets the real line in the interval $[-4,4]$ and the imaginary axis in $[-2,2]$. From David Wright.}

\medskip

Presumably the hyperbolic metric on ${\cal R}\subset \IC$ is equivalent to the Teichm\"uller metric.

\subsection{Discrete but not free points of ${\cal P}(F_2)$.}  Recently Akiyoshi, Oshika, Parker, Sakuma \& Yoshida,  \cite{AOPSY},  and also  Aimi, Lee, Sakai \& Sakuma \cite{ALSS}, have given complete proofs of the classification of all of these groups,  as announced by Agol, \cite{Agol}.  Historically, Riley's study of Kleinian groups generated by two parabolic transformations \cite{R1,R2},  and especially the construction of
the complete hyperbolic structure on the  figure-eight knot complement, inspired Thurston to establish the uniformisation theorem for  knots which are not either a torus knot or a satellite,  and more generally Haken manifolds, \cite{Thurston}. 

\medskip

The classification theorem is as follows - it basically shows that all Kleinian groups generated by two parabolics are near relatives to the two bridge knot and link groups.

\begin{theorem}  A non-free Kleinian group $\Gamma$  is generated by two  
parabolic elements if and only if one of the following holds.
\begin{itemize}
\item $\Gamma$ is conjugate to the hyperbolic 2-bridge knot or link group, 
\[ G(r)=\pi_1(\IS^3\setminus K(r)) \]
for some rational slope $r=q/p$, where $p$ and $q$ are relatively prime integers with $q\neq \pm1 \mbox{(mod $p$)}$
\item $\Gamma$ is conjugate to the Heckoid group, $G(r; n)$, for some rational $r$ and some $n=p/2$,  $p\geq 3$.
\end{itemize}
Further,  every hyperbolic 2-bridge link group $G(r)$ has precisely two parabolic generating pairs up to equivalence,  and 
every Heckoid group $G(r; n)$ has a unique parabolic generating pair up to
equivalence.
\end{theorem}

Thus in the complement of the Riley slice are isolated points accumulating at the boundary of the Riley slice.  Each point yields a representation $[\rho]\in {\cal P}(F_2)$ which is not free.  For instance,  the arithmetic points of the complement of the Riley slice are  the following two bridge knots and links,  \cite{GMcM}.
\begin{itemize}
\item $z_\rho=\frac{1+i\sqrt{3}}{2}$:   the figure of $8$ knot group, a two-bridge link of slope $\Big(\frac{5}{3}\Big)$.  It has index $12$ in $PSL(2,{\cal O}_{3})$.  It is the only arithmetic knot,  \cite{Reid}.
\item $z_\rho=1+i$:  the Whitehead link group, a two-bridge link of slope $\Big(\frac{8}{3}\Big)$and has index $12$ in the Picard group $PSL(2,{\cal O}_{1})$.
\item $z_\rho=\frac{3+i\sqrt{3}}{2}$:   the two-bridge link of slope $\Big(\frac{10}{3}\Big)$ and has index $24$ in $PSL(2,{\cal O}_{3})$. 
\item $z_\rho=\frac{1+i\sqrt{7}}{2}$:   the two-bridge link of slope $\Big(\frac{12}{5}\Big)$. 
\end{itemize}
Additionally there are the $(p,\infty,\infty)$-triangle groups,  $p\geq 3$ corresponding to the real points $z_\rho=-4\cos^2(\frac{\pi}{p})$. 
 
 \section{Main result.}

For the remainder of the paper we will be proving the following theorem concerning the complement of the Riley slice,  roughly it shows that `every group is everywhere'.  Generically (a dense $G_\delta$ set) a point in the Jordan domain  $\Omega=\IC\setminus {\cal R}$ represents a free group.  This occurs if,  for instance, $\gamma_\rho$ is not an algebraic integer. However we show below in Theorem \ref{fix} that a group $\langle f,g\rangle$ is not free if and only if $\gamma(f,g)$ is a fixed point of a word polynomial (defined below).

\begin{theorem} Let $\gamma_\rho\in \IC\setminus \{0\}$ and let $\Gamma_{\gamma_\rho}=\langle f,g\rangle$ where $f$ and $g$ are as (\ref{fg}).  Then exactly one of the following occurs.
\begin{enumerate} 
\item There is a neighbourhood $U$ of $\gamma_\rho$ such that for each $z\in U$ the group $\Gamma_z$ is geometrically finite,  discrete and free.
\item Let $z\neq 0$.  Then for every neighbourhood $V$ of  $\gamma_\rho$, there is $w\in V$ so that $\Gamma_z$ is a subgroup of $\Gamma_w$.
\end{enumerate}
\end{theorem}

\subsection{$\IZ_2$ extensions.} 
Before we move to the proof of this theorem we give another description of the Riley slice which is better suited to our proof.  Let 
$ \phi = \pm \left[\begin{array}{cc} 0 & 1/\sqrt{z} \\   \sqrt{z} & 0 \end{array}\right] $.
Then $\phi$ represents an elliptic element of order $2$ in $PSL(2,\IC)$ and with $f,g$ as at (\ref{fg}) (and dropping the subscript $\rho$) we have
\begin{equation}
g = \phi f \phi^{-1}, \;\;\;\; \gamma(f,\phi)=z,  \;\;\; z^2=\gamma(f,\phi f \phi^{-1} )=\gamma(f,g) = \gamma.
\end{equation}
The choice of square root here does not matter - only that a consistent choice is made. Notice then that $\langle f,g\rangle \subset \langle f,\phi\rangle$ and that the index of the first in the second is at most two.    Note also that the choice
$ \phi = \pm \left[\begin{array}{cc} 0 & i/\sqrt{z} \\ i  \sqrt{z} & 0 \end{array}\right] $ gives $\gamma(f,\phi)=-z$.

\begin{theorem}
$\overline{ {\cal R}} = \{\gamma(f,\phi)\in \IC \}$ where   $f$ is parabolic, $\phi$ is elliptic of order two and $\langle f,\phi \rangle$ is discrete and free on these generators.  
\end{theorem}
In fact ${\cal R} = \{\gamma(f,\phi)\in \IC \}$ as above with the additional assumption that the group is geometrically finite.

\begin{theorem}\label{psi} Let $h= \pm \left[\begin{array}{cc} a & b \\ c & d \end{array}\right]$, $ad-bc=1$ be an element of $PSL(2,\IC)$ and let $f$ be as at (\ref{fg}).  Then there is $\psi\in PSL(2,\IC)$,  elliptic of order two,  such that $\psi f \psi^{-1}=hfh^{-1}$ and $\gamma(f,h) = \gamma(f,\psi)$.
\end{theorem}
\noindent{\bf Proof.} We may suppose that $\gamma(f,h)=c^2\neq 0$ for otherwise the result is trivial.  We calculate with $\psi=  \pm \left[\begin{array}{cc} x & y \\ u & -x \end{array}\right]$, $-x^2-uy=1$,  that
\[ hfh^{-1}=\left(
\begin{array}{cc}
 1-a c & a^2 \\
 -c^2 & a c+1 \\
\end{array}
\right), \;\;\;\; \psi f \psi^{-1}=\left(
\begin{array}{cc}
 1-u x & x^2 \\
 -u^2 & u x+1 \\
\end{array}
\right) \]
whereupon the obvious choice is $x=a$ and $u=c$.  Then $y=-(1+a^2)/c$ gives us $\psi$. \hfill $\Box$

\medskip

This gives the following corollary (a more general result appears in \cite{GM1,MM}).

\begin{theorem}\label{thm6} Let $\rho\in R(F_2)$ be a representation with $\rho(a)$ parabolic.  Then there is $\tilde\rho\in{\cal R}(F_2)$ with $\tilde\rho(a)=\rho(a)$ and $\tilde\rho(b)$ elliptic of order two and $\gamma_\rho(a,b)=\gamma_{\tilde\rho}(a,b)$.  If $\rho$ is discrete,  then so is $\tilde\rho$.
\end{theorem}
\noindent{\bf Proof.} We may conjugate $\rho$ by $\alpha\in PSL(2,\IC)$ so as to assume $\rho^\alpha(a)=f$ is of the form (\ref{fg}) and put $\rho^\alpha(b)=h$.  We put $(\tilde\rho)^\alpha(b)=\psi$ as constructed in Theorem \ref{psi}. By conjugating back,  this resolves everything but the discreteness.  However,  if $\langle f,h \rangle$ is discrete,  then so is $\langle f, hfh^{-1}\rangle=  \langle f, \psi f\psi^{-1}\rangle$ and this last group is finite index in  $\langle f, \psi \rangle$.  The result follows. \hfill $\Box$

\subsection{Word polynomials and their semigroup.} One of the main tools we use to investigate the parabolic variety $P(F_2)$ and hence the Riley slice is the following remarkable property of commutators.
Let $m\geq 1$ and ${\bf s}=(s_i)_{i=1}^{m}\in (\IZ\setminus\{0\})^m$.  Define a word $w_{\bf s}\in F_2=\langle a,b\rangle$ by
\begin{equation}\label{ws}  w_{\bf s} = b a^{s_1} b^{-1} a^{s_2} b \cdots b^{\mp 1} a^{s_k} b^{\pm 1} \cdots a^{s_m} b^{\pm 1} \end{equation}
Thus the exponents of $a$ are given by ${\bf s}$ and the exponents of $b$ oscillate in sign -  it is immaterial whether the first exponent is positive or negative.

The next theorem was first proved in \cite{GM1} by a long induction,  but subsequently established within a more general framework in \cite{MM} via quaternion algebras and a generalisation of the polynomial Pell equation.   

\begin{theorem}[Word Polynomials]\label{words} Let ${\bf s}\in \IZ^m$.  Then there is a polynomial $P_{\bf s}\in \IZ[\gamma,\beta]$,  monic in $\gamma$, with the following property.  Let $\rho\in R(F_2)$ have principal character $\chi_\rho=(\gamma,\beta,\tilde{\beta})$.  
Then $\rho_w\in R(F_2)$ defined on $F_2\cong \langle u,v\rangle$ by 
\[\rho_w(u) = \rho(a),  \;\;\;\; \rho_w(v)=\rho(w_{\bf s}) \]
has principal character $  \big(P_{\bf s}(\gamma,\beta),\beta,\beta(\rho(w_{\bf s}))\big)$. Furthermore if $\rho$ is discrete,  then so are $\rho_w$ and  $\tilde{\rho}_w$.
 \end{theorem}
 
 Now with Theorem \ref{thm6} we can see how these polynomials act on the complex plane,  and in particular the Riley slice.  It is important to notice that if $\phi$ is elliptic of order two,  then {\em any} word   $w\in \langle f,\phi\rangle$ is a good word since the alternating sign condition is redundant, and since initial and terminal powers of $f$ are immaterial to the calculation of  $\gamma(f,w)$.  We therefore make the following definition. Suppose $w\in \langle a,b|b^2=1 \rangle$.  Write $w=a^{m_1}ba^{m_2}\cdots b a^{m_n}$, $m_i\neq 0$ for $2\leq i\leq n-1$, let ${\bf s}=(m_2,m_3,\ldots,m_{n-1})$ and define
 \[ p_w(z,\beta)=P_{{\bf s}}(z,\beta). \]
 
 \begin{theorem} Let $(\gamma,0,-4)$ be the principal character of an irreducible  representation $\rho\in R(F_2)$ and $w \in \rho(F_2)$,  then $(p_w(\gamma),0,-4)$ is the principal character of a representation $\rho_w\in R(F_2)$.  In addition $\rho_w$ is discrete if $\rho$ is.
 \end{theorem}
 \noindent{\bf Proof.} The group $\rho(F_2)$ is conjugate to the group $\langle f ,\phi \rangle$ with $f$ parabolic and $\phi$ of order $2$ and $\gamma(f,\phi)=\gamma)$ as above in \S 2.1 and so we may suppose the image is in fact this group.  If $w\in \langle f ,\phi \rangle$,  then as noted above, $\phi=\phi^{-1}$ the alternating sign condition in Theorem \ref{words} is redundant,  $w$  has polynomial $p_w$.  The group $\langle f,w\rangle$ is then discrete with $\langle f,\phi\rangle$ and has principal character $(p_w(\gamma),0,\beta(w))$.  Theorem \ref{thm6} then implies $(p_w(\gamma),0,-4)$ is a principal character,  discrete with $\langle f,w\rangle$. \hfill $\Box$
 
 \medskip
 
 \begin{corollary} If $p_w$ is a word polynomial and $z\in \IC\setminus\{0\}$ has $\langle f,g\rangle$ discrete as at (\ref{fg}),  then $p_w(z)$ is also a point of discreteness. \end{corollary}
 \noindent{\bf Proof.} Observe that $z=\gamma(f,\phi)$. \hfill $\Box$
 
 \medskip
 
 Let us give a few simple examples of these polynomials when $f$ is parabolic and $h$ is arbitrary.  Set $\gamma=\gamma(f,h)$. 
 
 \medskip

{\bf Table 1. Some examples of word polynomials. }\\
\begin{tabular}{|c|c|c|c|c|}
\hline
 Polynomial & word  & &    Polynomial & word   \\
 \hline
$\gamma^2$ & $ hfh^{-1} $   && 
$4 \gamma^2 $ & $hf^2h^{-1}$ \\
\hline
$  \gamma(1 -\gamma)^2 $ & $ hfh^{-1}fh$ && 
$\gamma (1+\gamma)^2$ & $ hfh^{-1}f^{-1} h$ \\
\hline 
$\gamma   (1-2 \gamma)^2$ & $ hfh^{-1}f^2h$ &&
$\gamma(1-4 \gamma)^2  $ & $ hf^2h^{-1} f^2h$\\
\hline
$  \gamma^2 (2 -\gamma )^2 $ & $hfh^{-1}fhfh^{-1}$ &&
$ \gamma ^2 \left(324 \gamma ^2-180 \gamma +25\right) $ & $hf^3h^{-1}f^{-3}hf^3h^{-1}$\\
\hline
\end{tabular}
 
\bigskip

There is a semigroup operation on the words $w_{\bf s}$ which is easiest explained as follows:
\begin{equation}\label{4}
{\bf s}*{\bf t} \leftrightarrow w_{\bf s}*w_{\bf t} = w_{\bf s} a^{t_1} w_{\bf s}^{-1} \cdots w_{\bf s} a^{t_m} w_{\bf s}^{-1}
\end{equation}
 That is we replace every instance of $b$ in $w_{\bf t}$ by $w_{\bf s}$.  Thus ${\bf s}*{\bf t}\in \IZ^{mn+m+n}$ when ${\bf s}\in \IZ^m$ and ${\bf t}\in \IZ^n$.

It is not too hard to see this corresponds to polynomial composition. In general $P_{{\bf s}*{\bf t}}(\gamma,\beta)=P_{\bf s}(P_{\bf t}(\gamma,\beta),\beta)$,  however we have $\beta=0$ and which simplifies things a bit to the following. 
 \begin{corollary} If $p_w$ is a word polynomial and $z\in \IC\setminus\{0\}$ has $\langle f,g\rangle$ discrete as at (\ref{fg}),  then $p_w(z)$ is also a point of discreteness. \end{corollary}
\begin{theorem}\label{thm9} If $p_w$ and $p_v$ are  word polynomials,  then so is $p_w\circ p_v$.  In fact $p_w\circ p_v=p_{w*v}$.
\end{theorem}

\subsection{The Julia set of the word semigroup is the exterior of the interior of the Riley slice.}
Let 
\[ {\cal G}=\{p_w(z) : \mbox{$w$ is a good word}\}\]
We have seen that ${\cal G}$ is a family of polynomials closed under composition  - thus forming a semigroup. It turns out that ${\cal G}$ is infinitely generated as there are infinitely many polynomials of prime order.  We wish to use the mixing properties of the semigroup ${\cal G}$ on its Julia set to establish our main result.  The dynamical theory of semigroups of polynomial mappings,  or more generally semigroups of rational mappings,  was initiated in \cite{HM1} as a generalisation of the theory of iteration of rational mappings.  Since then it has grown more broadly,  see for instance the work of Sumi, Stankewicz and Urba\'nski \cite{SS,S,SU} etc. We begin with a few definitions.  

\medskip

The {\em Fatou set} of a polynomial semigroup $G$,  $F(G)$,  is the set of all $z\in \IC$  which admit a neighbourhood  $z\in U$ so that the family $G|U$ is normal on $U$ -- every sequence contains a locally uniformly convergent subsequence.  The {\em Julia set} of $G$ is $J(G)=\IC\setminus F(G)$.  By definition $F(G)$ is open and hence $J(G)$ is closed. Unlike the case of polynomial iteration these sets are not completely invariant but we summarise some basic facts in the next lemma.

\begin{lemma}\label{genfacts} Let $G$ be a rational semigroup.
\begin{itemize}
\item  For any $g\in G$,
\begin{itemize}
\item $g(F(G))\subset F(G),$ and $g^{-1}(J(G))\subset J(G)$
\item $F(G)\subset F(\langle g\rangle)$,  and $J(\langle  g \rangle)\subset J(G)$
\end{itemize} 
\item If $J(G)$ contains at least three points, then $J(G)$ is a perfect set.
\item  If there is some $g\in G$ such that $\deg(g)\geq 2$  and $J(G)$ contains at least three
points, then $J(G)$ is the smallest closed backward invariant set containing
at least three points. Here we say that a set $A$ is backward invariant under
$G$ if for each $g\in G,$ $g^{-1}(A)\subset A$.
\item If $J(G)$ contains at least three points, then
\[ J(G)=\overline{\{z\in \oC:   \mbox{$z$ is a repelling fixed point of some $g\in G$}\}}\]
\end{itemize}
 \end{lemma}

The next result from \cite{HM2} will be important to us.  It is a consequence of the fact that a finitely generated polynomial semigroup has uniformly perfect Julia set.  It gives a criteria for the Julia set to have nonempty interior.

\begin{lemma}\label{open} Let $p(z)$ and $q(z)$ be polynomials of degree at least two.  Suppose that $z_0$ is a common fixed point of $p$ and $q$,  $p(z_0)=q(z_0)=z_0$,   that  $p$ has $z_0$ as a super-attracting fixed point - $p'(z_0)=0$, and that $z_0\in J(\langle q \rangle)$.  Then $J(\langle p,q \rangle)$ contains a neighbourhood of $z_0$.
\end{lemma}

Now applying these results in our setting we have a first result. 

\begin{corollary}\label{fij} The Julia set of the word semigroup ${\cal G}$ contains a neighbourhood $V_{\cal G}$ of $0$.
\end{corollary}
\noindent{\bf Proof.} ${\cal G}$ contains the two polynomials $z\stackrel{p}{\mapsto } z^2$ and $z\stackrel{q}{\mapsto } z(1-z)^2$ both with nontrivial Julia sets.  The origin is a super-attracting fixed point of  $p$.  We compute that  $q'(0)=1$ and therefore $0$ is a parabolic fixed point lying in the Julia set of $q$,  illustrated below in Figure 2.  The claim now follows from Lemma \ref{open}. \hfill $\Box$

 \scalebox{0.45}{\includegraphics[viewport= -5 360 580 790]{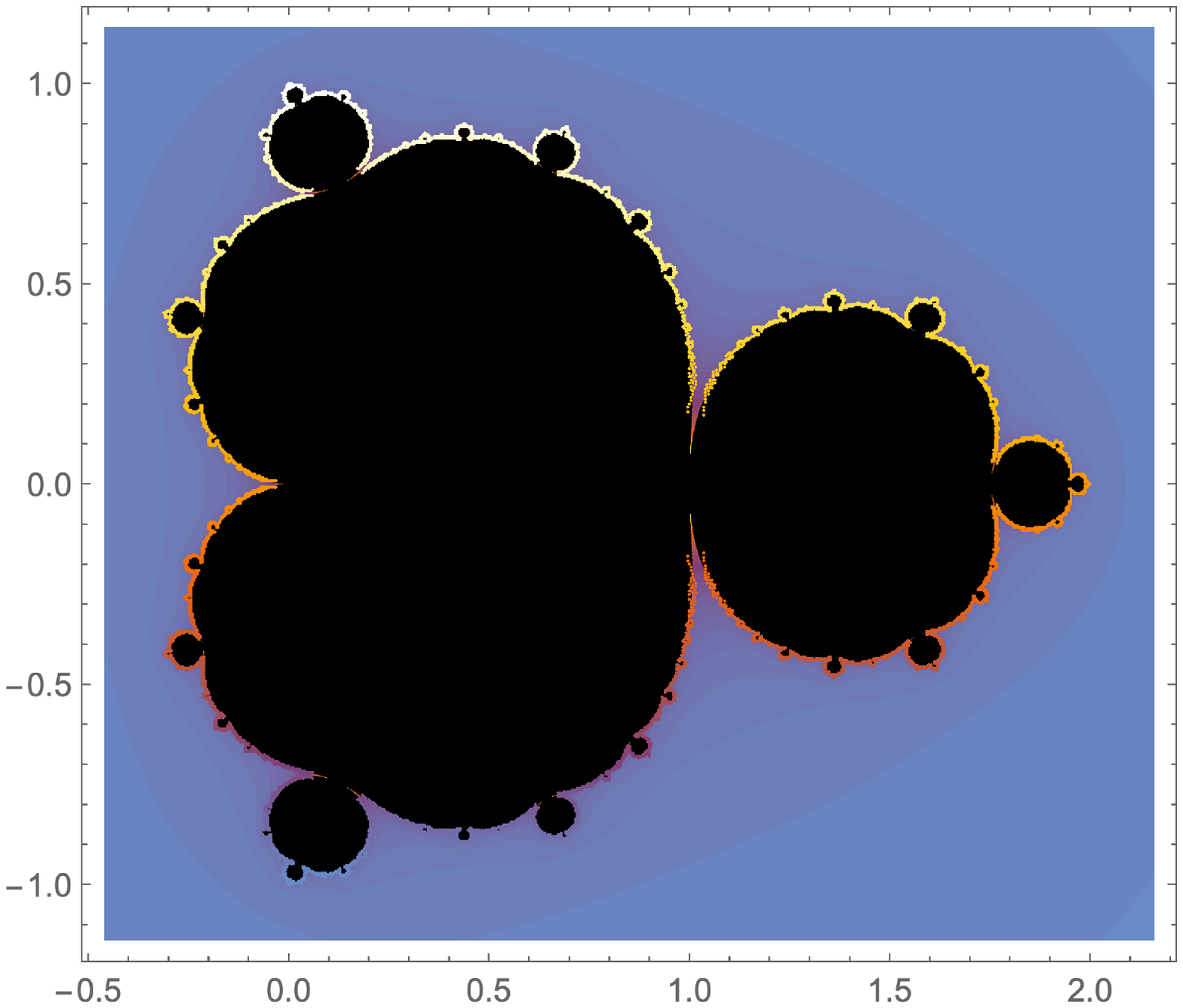}}\\
{\bf Figure 2.}  {\em The filled Julia set of the polynomial $\gamma\mapsto \gamma  (1-\gamma )^2$.  This set lies in the complement of the Riley slice (Lemma \ref{fj}) and has the $2$-congruence subgroup of the modular group ($\gamma(f,\phi)=z_\rho=2$) as a periodic cycle of period $1$,   the figure-eight knot group ($z_\rho=\frac{1+i\sqrt{3}}{2}$)  as a periodic cycle of period $2$ and both the Gaussian integer group ($z_\rho=i$) and the modular group ($z_\rho=1$) a eventually periodic points of period $1$.  See \S 3 for a more general discussion of these periodic points. }

\bigskip
Next we show that no polynomial in ${\cal G}$ has a root in the Riley slice.

\begin{theorem} \label{roots} Let $p\in {\cal G}$ and $z_0\in {\cal R}$, the Riley slice.  Then $|p(z_0)| \geq 1$. 
\end{theorem}
\noindent{\bf Proof.} Suppose that $|p(z_0)|<1$. First,  as $z_0\in{\cal R}$,  $z_0\neq 0$ and $(z_0,0,-4)$ is the principal character of a discrete representation whose image is free on the two generators.  Then there is a neighbourhood $U_0$ of $z_0$ so that $\overline{p(U_0)}\subset \ID(0,1)$.  There is $0\neq z_1\in U_0\cap {\cal R}$ so that $|p(z_1)|<1$.  Then $(p(z_1),0,-4)$ is the principal character of a discrete representation $\rho\in R(F_2)$.  This group is non-elementary as it contains a parabolic whose fixed point is not one of those of the element of order two as $\gamma_\rho\neq0$. Hence the group contains two parabolics with distinct fixed points.  However the Shimitzu-Leutbecher inequality tells is that for a nonelementary group with containing a parabolic element $f$, for each $g$ we have $|\gamma(f,g)|\geq 1$ unless $\gamma(f,g)=0$. \hfill $\Box$

\medskip

\begin{corollary} \label{Fatou} $F(\cal G) \supset {\rm int}({\cal R})$.
\end{corollary}
\noindent{\bf Proof.} If $z_0\in {\rm int}({\cal R})$,  then there is a neighbourhood $V_0$ of $z_0$ containing only discrete nonelementary representations. Each $p\in {\cal G}$ must have $p(V_0)\subset \IC\setminus \ID$, that is the family must omit the unit disk on $V_0$.  This is a standard normality criterion for families of analytic functions (c.f. Montel's Theorem). \hfill $\Box$

\medskip
\noindent{\bf Remark.} Actually the Shimitzu-Leutbecher inequality is a simple consequence of the fact that $z\mapsto z^2 \in {\cal G}$.  If $|z_0|<1$,  then 
\[ z_1 = \gamma(f,\phi f \phi)=z_0^2,  z_2=\gamma(f,\phi f \phi^{-1} f \phi f^{-1} \phi^{-1}) =  \gamma(f,\phi f \phi)^2 = z_0^4, \ldots \]
and iterating this procedure provides a sequence of traces converging $2$.  This is easily seen to contradict discreteness.

\medskip

We begin to see that the word semigroup is quite remarkable - it is infinitely generated, closed under composition,  and all the roots of all the elements lie in a bounded region of $\IC$.

\begin{theorem}\label{JS} The Julia set of  ${\cal G}$ is the exterior of the interior of the Riley slice.
\end{theorem}
\noindent{\bf Proof.} Since the Julia set is closed,  since the discrete representations outside the Riley slice are isolated, and having in hand Corollary \ref{Fatou}, it is enough to show that for any $(z_0,0,-4)$ which is the principal character of a non-discrete representation ${\cal G}$ is not  normal in a neighbourhood of $z_0$.  Then suppose that $\langle f,\phi \rangle$ is not discrete, $f$ parabolic, $\phi$ of order two, and $z_0=\gamma(f,\phi)$.  Since $\langle f,\phi \rangle$ is not discrete there is a sequence $w_i\in \langle f,\phi \rangle$ with $w_i \to identity$ in $SL(2,\IC)$.  Then
\[ p_{w_i}(z_0) = \gamma(f,w_i) \to 0\]
 Then for all $i$ sufficiently large, $p_{w_i}(z_0) \in V_{\cal G}$,  where $V_{\cal G}$ is given in Corollary \ref{fij}.  It follows that there is $i_0$ and an open $U_0$ about $z_0$ so that $p_{w_{i_0}}(U_0)\subset V_{\cal G}$.  Since ${\cal G}$ is not normal on any open subset of $V_{\cal G}$ and since ${\cal G}\circ p_{w_{i_0}}\subset {\cal G}$,  the semigroup ${\cal G}$ cannot be normal on $U_0$ or on any open subset of it. Thus $z_0\not\in F({\cal G})$ and hence $z_0\in J({\cal G})$ and this completes the proof. \hfill $\Box$

\medskip

We believe this is the first example of a polynomial semigroup whose Julia set is a Jordan domain in $\IC$.  Is it a quasi-circle ?

\subsection{Completion of the proof.}

Let $z_0\in \IC\setminus \{0\}$ and let $\lambda \in \overline{ \IC\setminus {\cal R}}$, $\lambda \neq 0$, a point in the complement of the interior of the Riley slice.  Our target group is
\[ \Gamma_0= \Big\langle \left[\begin{array}{cc} 1 & 1 \\ 0& 1 \end{array}\right],\left[\begin{array}{cc} 1 & 0 \\ z_0  & 1 \end{array}\right] \Big\rangle \]
whose principal character is $(z_0^2,0,0)$,  and we are seeking a supergroup of $\Gamma_0$ near 
\[ \Gamma_\lambda= \Big\langle \left[\begin{array}{cc} 1 & 1 \\ 0& 1 \end{array}\right],\left[\begin{array}{cc} 1 & 0 \\ \lambda  & 1 \end{array}\right] \Big\rangle . \]

Let $\langle f,\phi_\lambda \rangle\subset PSL(2,\IC)$ be the $\IZ_2$ extension of $\Gamma_\lambda$ with principal character $(\lambda,0,-4)$, so
\[ \phi_\lambda f\phi_\lambda ^{-1}= \left[\begin{array}{cc} 1 & 0 \\ \lambda  & 1 \end{array}\right] \]

The closure of the backward orbit of $z_0$ under the polynomial semigroup ${\cal G}$,
\[ \overline{{\cal G}^{-1}(z_0)}= \overline{\{p_w^{-1}(z_0):p_w\in {\cal G}\}} \]
is closed, contains more than three points (use the identified polynomials in Table 1.) and hence contains the Julia set of ${\cal G}$ by Lemma \ref{genfacts}. This Julia set is the closure of the complement of the Riley slice by Theorem \ref{JS}. In particular every neighbourhood of $\lambda$ contains such a preimage. 

\medskip

We have now shown that every neighbourhood of $\lambda$ contains a point $\zeta$ for which there is a word $w$ in the abstract group $\IZ*\IZ_2$ for which $p_w(\zeta)=z_0$ for any representation for which the image of the $\IZ$ factor is represented by a parabolic (and the image of $\IZ_2$ is not trivial).   Thus let $\langle f,\phi_\zeta \rangle \subset PSL(2,\IC)$ have principal character $(\zeta,0,-4)$.  The subgroup $\langle f,w_\zeta \rangle$,  where we replace the abstract word with instances of $f$ and $\phi_\zeta$, has principal character 
\[ (p_w(\gamma(f,\phi_\zeta)),0,\beta(w_\zeta))=(p_w(\zeta),0,\beta(w_\zeta))=(z_0,0,\beta(w_\zeta)) \]
Therefore the subgroup $\langle f,w_\zeta f w_\zeta^{-1}\rangle$ has principal character $(z_0^2,0,0)$ and so is conjugate to $\Gamma_0$.  Finally we observe that the word $w_\zeta f w_\zeta^{-1}$ has an even number of occurrences of $\phi_\zeta (=\phi_{\zeta}^{-1})$.  Thus $w_\zeta f w_\zeta^{-1} \in \langle f,\phi_\zeta f \phi_{\zeta}^{-1} \rangle$ and hence
\[ \langle f,w_\zeta f w_\zeta^{-1}\rangle \subset \langle f,\phi_\zeta f \phi_{\zeta}^{-1} \rangle \]
On the left-hand side here we have a conjugate of $\Gamma_0$ and on the right-hand side we have a group generated by two parabolics arising from the point $\zeta$.

\medskip

This completes the proof. \hfill $\Box$

\bigskip

We have the easy consequence in a very special case.

\begin{corollary} Let $\lambda$ be a point in the complement of the Riley slice.  Then every neighbourhood of $\lambda$ contains a point $\zeta$ so that $\Gamma_\zeta$ contains the figure-eight knot group as an infinite index subgroup.
\end{corollary}

There have been attempts to describe the Julia sets of polynomial semigroups algorithmically,  see \cite{SS2}. In our case we have an alternative description of the Julia set and so can make some comparisons.  We lexicographically ordered the first few thousand words,  found the polynomials and the calculated their roots,  and also the preimages of the figure-eight knot group.  These are shown below in Figure 3. In each case we obtained around 20-40 thousand distinct points.  There appears to be a strong clustering to the real axis,  and it is very difficult to get near the boundary with these ideas.

\scalebox{0.5}{\includegraphics[viewport= 30 0 380 300]{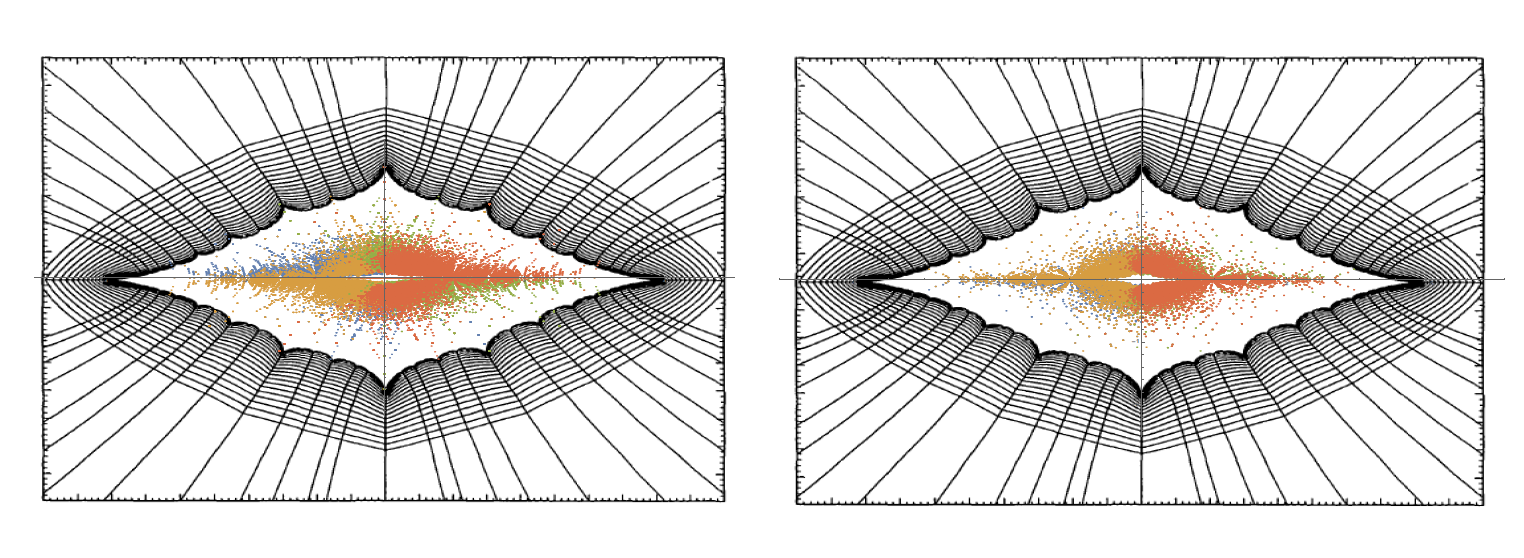}}\\
{\bf Figure 3.}  {\em The Riley slice. Left with polynomial roots.  Right with backward orbit of Figure-eight knot group and its symmetries}

\section{Iteration and Nielsen classes.}

In \cite{ALSS,AOPSY} it is proved,  as a consequence of their classification,  that a discrete group generated by two parabolics has either one or two Nielsen classes of parabolic generating pairs.  This actually gives a new criterion for discreteness. Here we wish to show that this is far from the case for a dense family of groups in the complement of the Riley slice.

\begin{theorem}\label{NC} Let $N\geq 1$ be an integer.  Then the set of points in the complement of the Riley slice representing a group generated by two parabolics with at least $N$ distinct Nielsen classes of parabolic generating pairs is dense.
\end{theorem}

Let $w\in \langle f,g: f^\infty=g^2 = \cdots =1\rangle $,  and $*$ the semigroup operation as described at (\ref{4}).  Then  $\gamma(f,w)=p_w(\gamma(f,g))$ for a monic polynomial $p_w$ with integer coefficients.   We define the $n$-fold product,
\[ *^nw=w*w*\cdots*w, \hskip10pt \mbox{$w$ occurs $n$ times}, \]
and observe that
\[ p_{*^nw}(z) = p^{(n)}(z)=p(p(\cdots p(z)\cdots)) \]
is the $n^{th}$ iterate of $p$.  Let $\gamma_1,\gamma_2,\ldots,\gamma_{2N-1},\gamma_{2N}$ be a periodic cycle of length $2N$.  Thus $z_i=p_{*^i}(z_1)$, $z_i\neq z_j$ if $i\neq j$,  $p_w(z_{2N})=z_1$ and there are at least $N$ distinct values in the set \{ $\gamma_i^2:i=1,\ldots 2N\}$.   Here the cycle $\{z_i\}$ may be in the Julia set or in the bounded components of the Fatou set of the poynomial $p_w$.  We first observe the following easy lemma,  see Figure 2,  which is also a corollary of Theorem \ref{roots}.

\begin{lemma}\label{fj} Let $w$ be a good word with polynomial $p_w$.  Then the filled Julia set of $p_w$,  the complement of the attracting basin of the super-attracting point  $\infty$, lies in the complement of the Riley slice.
\end{lemma} 
\noindent {\bf Proof.} The boundary of the filled Julia set is the Julia set.  The semigroup ${\cal G}$ contains $p_w$ and cannot be normal on its Julia set.  Thus every Julia set lies in the complement of the Riley slice and since this complement is simply connected the result follows. \hfill $\Box$

\bigskip

We say that the group $\Gamma_\gamma$,
\[ \Gamma_\gamma = \Big\langle  \left[\begin{array}{cc} 1 & 1 \\ 0& 1 \end{array}\right],   \left[\begin{array}{cc} 1 & 0 \\ {\gamma} & 1 \end{array}\right]\Big\rangle,    \]
is periodic of period $m$,  or eventually periodic of period $m$, if there is a good word $w$ and $\gamma$ is such a point for $p_w$.

\begin{theorem} Let $\Gamma_\gamma$ be $2N$ periodic.  Then $ \Gamma_\gamma$ has at least $N$ Nielsen distinct pairs of parabolic generators.
\end{theorem}
\noindent{\bf Proof.}   Let $w$ be a word and suppose $\gamma$ lies on a periodic cycle of length $2N$ of $p_w$, Denote this cycle by $\gamma=\gamma_1\neq  \ldots \neq \gamma_{2N}$,  $p_w(\gamma_{2N})=\gamma$.  Let $\Gamma^{(2)}=\langle f,\phi \rangle$ be the $\IZ_2$-extension of $\Gamma_\gamma$ with principal character $(\gamma,0,-4)$,  as in \S 2.1.  Then consider the subgroups $\langle f,*^i w\rangle$.  They have principal character $(\gamma_i,0,\beta(*^iw))$.  For $i\geq 1$ the groups
\[ \Gamma_i =  \langle f,( *^iw )\circ f\circ  (*^i w)^{-1} \rangle \subset \Gamma_\gamma\]
as there is an even number of occurrences of $\phi$ in the word $( *^iw )\circ f\circ  (*^i w)^{-1}$.  The principal character of these groups is $(\gamma_i^2,0,0)=(p_{w}^{(i)}(\gamma)^2,0,0)$ and since there are at least $N$ distinct values here,  there are at least $N$ distinct Nielsen classes as $\gamma$ is an invariant of the Nielsen class.  Now notice that $\gamma_{2N+1}=\gamma$ and hence $\Gamma_{2N+1} \subset \Gamma_\gamma$ is actually conjugate to $\Gamma_\gamma$.  Of course this conjugacy is a bijective homomorphism.  It follows that the chain of iterates and subgroups is a chain of equalities,
\[ \Gamma_\gamma = \Gamma_{\gamma_1} =  \Gamma_{\gamma_2} =\cdots = \Gamma_{\gamma_{2N}} = \Gamma_\gamma .\]
We have identified at least $N$ distinct Nielsen classes of generators and so the theorem is proved. \hfill $\Box$
 
\bigskip

Since high order periodic points are dense in the Julia set of a polynomial,  and since Julia sets of word polynomials are dense in the complement of the Riley slice,  Theorem \ref{NC} follows directly.

 \bigskip
 
 Finally we establish the following related theorem alluded to earlier.
 
 \begin{theorem}\label{fix} A group  $\Gamma_\gamma$ generated by two parabolics is not free if and only if $\gamma$ is a fixed point of a word polynomial.  Further,  if $\gamma$ is eventually periodic,  then $\Gamma_\gamma$ is not free.
 \end{theorem}
 \noindent {\bf Proof.} If $\gamma$ is eventually periodic for, say $p_w$,  then we may iterate until we get to some $n$ with $p_{w}^{\circ n}(\gamma)=\gamma_0$,  which is a fixed point for this polynomial.  Since subgroups of free groups are free,  it is enough to show that $\Gamma_{\gamma_0}$ is not free where $\gamma_0$ is fixed by some word polynomial $p_{w_0}$.  As above we find that the group generated by two parabolics $\langle f,g \rangle$, $\gamma_0=\gamma(f,g)$ is conjugate to the group generated by two parabolics $\langle f,w_0f w_{0}^{-1} \rangle$,  $\gamma(f,w_0f w_{0}^{-1})=p_{w_0}(\gamma_0)=\gamma_0$,  and so this latter group is both a subgroup and a conjugate of the former.  Thus $\langle f,w_0f^{-1}w_0\rangle=\langle f,g\rangle$.  In particular $w_0\in \langle f,w_0f^{-1}w_0\rangle$.  In these generators $w_0$ can be expressed as a word containing an even number of $w_0$'s.  It follows the group is not free. 
 Next,  suppose that  $\Gamma_\gamma$ is not free.  Then the $\IZ_2$ extension $\langle f,\phi\rangle$, $\gamma=\gamma(f,\phi)$,  is not free on its generators either.  Let $u\in \langle f,\phi\rangle$ be a maximally reduced  nontrivial word representing the identify.  We can write $u$ in the form
 \[ u = f^{a}w f^{b}\phi f^{c} = identity\]
 where $w$  starts and ends in $\phi$ and possibly $a$ and/or $c$ are $0$.  Now $w$ is a good word with polynomial $p_w$, but also as words in this particular group, $w= f^{-a-c}\phi f^{-b}$.  We calculate that
 
\[ p_w(\gamma)= \gamma(f,w)=\gamma(f,f^{-a-c}\phi f^{-b})=\gamma(f, \phi )=\gamma. \]
This is what we wanted to prove.
 \hfill $\Box$

\bigskip 
 \noindent G. Martin, Massey University,   New Zealand,  g.j.martin@massey.ac.nz

\begin{thebibliography}{99}
\bibitem{Agol} I. Agol, {\em The classification of non-free 2-parabolic generator Kleinian groups}, Slides of talks given
at Austin AMS Meeting and Budapest Bolyai conference, July 2002, Budapest, Hungary.
\bibitem{ALSS} S. Aimi, D. Lee, S. Sakai, and M. Sakuma, {\em Classification of parabolic generating pairs of Kleinian
groups with two parabolic generators}, preprint.
\bibitem{AOPSY} H. Akiyoshi, K. Ohshika, J. Parker, M. Sakuma and H. Yoshida, {\em Classification of non-free
Kleinain groups generated by two parabolic transformations}, preprint.
\bibitem{ASWY} H. Akiyoshi, M. Sakuma, M. Wada and Y. Yamashita, {\em  Punctured torus groups and two bridge knot groups (I)},  Lecture Notes in Mathematics 1909, Springer-Verlag Berlin Heidelberg, 2007.
\bibitem{B} A. Beardon, {\em The geometry of discrete groups},
Springer--Verlag, 1983.  
\bibitem{Bers}  L. Bers, {\em Uniformization, moduli and Kleinian groups}, Bull. London Math. Soc., {\bf 4}, (1972), 257--300.
\bibitem{FK} R. Fricke and F. Klein, {\em Vorlesungen \"uber die Theorie der automorphen Functionen}, Chapter 2, Teubner, Leipzig, 1897. 
\bibitem{GM1} F. W. Gehring and G. J. Martin {\em Commutators, collars and the
geometry of M\"{o}bius groups}, J. d'Analyse  Math. {\bf 63} (1994) 175 -- 219.
\bibitem{GMcM} F. W. Gehring, C. Maclachlan and G. J. Martin {\em Two-generator
arithmetic Kleinian groups II},  Bull. London Math. Soc., {\bf 30}, (1998), 258
-- 266.
\bibitem{Gill} J. Gilman, {\em The structure of two-parabolic space: parabolic dust and iteration}, Geom. Dedicata,
{\bf 131}, (2008), 27--48.
\bibitem{HM1} A. Hinkkanen and G. J. Martin, {\em The dynamics of semigroups of rational functions I}, Proc. London Math. Soc., {\bf 73}, (1996),  58--84,. 
\bibitem{HM2} A. Hinkkanen and G.J. Martin, {\em Julia sets of rational semigroups},  Math. Z., {\bf 222}, (1996), 161--169.
\bibitem{Jorgensen} T. J\o rgensen, {\em On discrete groups of
M\"obius transformations}, Amer. J. Math.,
{\bf  98}, (1976), 739--749.
\bibitem{KS} L. Keen and C. Series, The Riley slice of Schottky space, Proc. London Math. Soc. 69 (1994), 72--90.
\bibitem{KoS}  Y. Komori and C. Series, The Riley slice revisited, The Epstein birthday schrift, 303--316,
Geom. Topol. Monogr., 1, Geom. Topol. Publ., Coventry, 1998.
\bibitem{Kra}  I. Kra, {\em Horocyclic coordinates for Riemann surfaces and moduli spaces I: Teichmiiller and
Riemann spaces of Kleinian groups}, Amer. J. Math., {\bf 3}, (1990), 499--578.
\bibitem{LU}  R. C. Lyndon  and J. L. Ullman, {\em  Groups Generated by two Parabolic Linear Fractional Transformations}, Canadian J. Math., {\bf 21},  (1969), 1388--1403.
\bibitem{MM} T. H. Marshall and G. J. Martin {\em Polynomial Trace Identities in $SL(2,{\bf C})$,  Quaternion Algebras, and Two-generator Kleinian Groups},  arXiv 1911.11643 (2019).
\bibitem{OM} K. Ohshika and H. Miyachi, {\em Uniform models for the closure of the Riley slice}, In the tradition
of Ahlfors-Bers. V, 249 -- 306, Contemp. Math., {\bf 510}, Amer. Math. Soc., Providence, RI, 2010.
\bibitem{Reid} A. W. Reid {\em Arithmeticity of knot complements},  J. London Math. Soc., {\bf 43},  (1991),   171--184.
\bibitem{R1} R. Riley, {\em A quadratic parabolic group}, Math. Proc. Cambridge Philos. Soc. 77 (1975), 281--288.
\bibitem{R2} R. Riley, {\em Parabolic representations of knot groups. II}, Proc. London Math. Soc. 31 (1975),
495--512.
\bibitem{Scott} G.P. Scott, {\em Finitely generated 3-manifold groups are finitely presented},  J. London Math. Soc.,   {\bf 6} (1973), 437–440.
\bibitem{SS}  R. Stankewitz and H. Sumi,  {\em Dynamical properties and structure of Julia sets of postcritically bounded polynomial semigroups}, Trans. Amer. Math. Soc., {\bf  363}, (2011),  5293--5319.
\bibitem{SS2}  R. Stankewitz and H. Sumi,  {\em Random backward iteration algorithm for Julia sets of rational semigroups},  Discrete Contin. Dyn. Syst., {\bf  35},  (2015),  2165--2175.
\bibitem{S}  H. Sumi, {\em Skew product maps related to finitely generated rational semigroups},  Nonlinearity, {\bf 13} (2000),  995--1019. 
\bibitem{SU}  S. Hiroki and M. Urba\'nski,   {\em Measures and dimensions of Julia sets of semi-hyperbolic rational semigroups},  Discrete Contin. Dyn. Syst., {\bf  30}, (2011),  313--363.
\bibitem{Thurston}  W. Thurston, The geometry and topology of three-manifolds, available from
http://library.msri.org/books/gt3m/
\end{thebibliography}
\end{document}